\documentclass[12pt,reqno]{amsart}
\usepackage{amssymb}

\usepackage{amsmath,latexsym,amsfonts}
\usepackage{graphicx}
\usepackage{srcltx}
\usepackage{amscd}
\usepackage{hyperref}

\newtheorem{lemma}{Lemma}[section]
\newtheorem{coro}[lemma]{Corollary}
\newtheorem{prop}[lemma]{Proposition}
\newtheorem{thm}[lemma]{Theorem}
\newtheorem{defn}[lemma]{Definition}
\newtheorem{rem}[lemma]{Remark}
\newtheorem{ex}[lemma]{Example}
\makeatletter\@addtoreset{equation}{section}
\renewcommand\theequation{\thesection.\@arabic\c@equation}
\newcommand{\func}{\operatorname}
\begin{document}

\title[Formulas for transversely K\"ahler foliations]
{Basic Dolbeault cohomology and  Weitzenb\"ock formulas on transversely K\"ahler foliations}
\author[S.~D.~Jung]{Seoung Dal Jung}
\address{Department of Mathematics\\
Jeju National University \\
Jeju 690-756 \\
Republic of Korea}
\email[S.~D.~Jung]{sdjung@jejunu.ac.kr}
\author[K.~Richardson]{Ken Richardson}
\address{Department of Mathematics \\
Texas Christian University \\
Fort Worth, Texas 76129, USA}
\email[K.~Richardson]{k.richardson@tcu.edu}

\begin{abstract}
We study basic Dolbeault cohomology  and find new 
Weitzenb\"ock formulas on a transversely K\"ahler foliation.  
We investigate conditions on mean curvature and Ricci curvature that impose
restrictions on basic Dolbeault cohomology.
For example, we prove that on a transversely K\"ahler foliation with positive transversal Ricci curvature, 
there are no nonzero basic-harmonic forms of type $(r,0)$, among other results.
\end{abstract}

\maketitle

 \renewcommand{\thefootnote}{}
\footnote {This paper was supported by  the National Research Foundation of Korea (NRF) grant funded by the Korea government (MSIP) (NRF-2015R1A2A2A01003491).
This work was also supported by a grant from the Simons Foundation (Grant Number
245818 to Ken Richardson).
}
\renewcommand{\thefootnote}{\arabic{footnote}}
\setcounter{footnote}{0}

\section{Introduction}
Let $\mathcal{F}$ be a foliation of a closed, smooth manifold $M$. One set of smooth invariants of $\mathcal F$ is the basic cohomology.
The basic forms of $(M,\mathcal F)$ are locally forms on the leaf space; that is, forms $\phi$ satisfying $X\lrcorner\phi = X\lrcorner d\phi =0$ for any vector $X$ tangent to the leaves, where $X\lrcorner$ denotes the interior product with $X$.
Basic forms are preserved by the exterior derivative and are used to define basic de Rham cohomology groups $H_B^*$. In 
general these groups can be infinite-dimensional.
A Riemannian foliation is a foliation on a smooth manifold such that the normal bundle
$Q=TM/T\mathcal F$
 may be endowed with a metric whose Lie derivative is zero along leaf directions. 
 We can always extend this metric to a metric $g_M$ on  $M$ that is called bundle-like, meaning that
the leaves of $\mathcal F$ are locally equidistant.  
From this additional structure on the foliation,
the basic Laplacian $\Delta_B$ is a version of the Laplace operator that preserves the basic forms. Many researchers have studied basic forms and the basic Laplacian on Riemannian foliations.
 It is well-known (\cite{Al}, \cite{EKH}, \cite{KT3},\cite{PaRi}) that on a closed oriented manifold $M$ with a transversally oriented Riemannian foliation $\mathcal F$, $H_B^r\cong \mathcal H_B^r$, where $\mathcal H_B^r={\rm ker}\Delta_B$ is finite dimensional.  
 Because of this Hodge theorem, researchers have been able to show relationships between curvature bounds and basic cohomology.
 In \cite{Heb}
 J. Hebda proved that a lower bound on transversal Ricci curvature for a Riemannian foliation of a compact manifold causes the space of leaves to be compact and the first basic cohomology group to be trivial.
 Another example relating the geometry and the topology came in 1991, when M. Min-Oo et al. \cite{MO} proved that if the transversal  curvature operator of $(M, \mathcal F)$ is positive-definite, then the cohomology $H_B^r=\{0\} \ (0<r<q)$;  that is, any basic-harmonic $r$-form is trivial. There are many other examples of known relationships between transversal curvature and basic cohomology.
 
We are interested in foliations that admit a transverse K\"ahler structure.
 These structures appear naturally in many situations. 
 For example, Vaisman manifolds, Calabi-Eckmann manfolds and Sasaki manifolds are not K\"ahler but carry a 
 transverse K\"ahler structure. In particular, the characteristic foliation on a Sasaki manifold is 
 transversely K\"ahler \cite[Section 2]{BGN} and is automatically homologically orientable \cite[Proposition 7.2.3]{BG},
 causing many facts about K\"ahler manifolds to carry over. Being homologically orientable
 implies that the foliation is taut; there exists a metric such that the leaves of the foliation are 
 minimal submanifolds. See \cite{Wo} for a treatment.
 For instance, in \cite{CMNY}, the authors were able  to 
 prove a hard Lefschetz theorem for compact Sasaki manifolds.
 In \cite{ChdLM}, the authors proved a similar result for cosymplectic manifolds.
 A major step in understanding the homologically orientable (taut) case was
 by A. El Kacimi-Alaoui in 
 \cite[Section 3.4]{EK}, where the standard facts about K\"ahler manifolds
 and their cohomology are
 shown to carry over to the transverse structures and basic cohomology.
 Also, L. A. Cordero and R. A. Wolak \cite{CW} studied the properties of the basic cohomology on a taut transversely K\"ahler foliation by using the differential operator $\Delta_T$, which is different from $\Delta_B$. In fact, $\mathcal F$ is minimal (taut) if and only if $\Delta_T =\Delta_B$. 
 Many known facts about K\"ahler forms $\omega$ carry over to the transversal geometry setting, but not all; for example, 
 if the foliation is not taut, the powers $\omega^k$ need not generate nonzero basic cohomology classes. 
In the general case, in \cite{JP} S. D. Jung and J. S. Pak proved a Weitzenb\"ock 
formula for the Dolbeault Laplacian and deduced among other results that when the mean curvature of the foliation is 
basic-harmonic, nonnegativity of the transverse Ricci curvature (and positivity at one point) implies that the first basic cohomology group is trivial.
In \cite{JP2}, they showed that under some positivity assumptions on transverse curvature on a complete transversely K\"ahler foliation, the second basic cohomology group is generated by the K\"ahler form.
In \cite{JL}, S. D. Jung and H. Liu proved that if the transversal Ricci curvature of a compact transversely K\"ahler 
foliation is nonpositive and negative at one point, then there are no nonzero basic transverse
holomorphic vector fields. G. Habib and 
L. Vezzoni showed a similar result by different methods in \cite[Theorem 3.3]{HV}. They used a 
twisted version of the basic Dolbeault Laplacian and calculated the corresponding Weitzenb\"ock formula.
In another recent paper \cite{KLW},
Y. Kordyukov, M. Lejmi, and  P. Weber showed that 
for a manifold with bundle-like, transversely K\"ahler metric for a taut foliation,
under some conditions on the basic cohomology of one-forms, 
the basic version of the Seiberg-Witten invariants must be nontrivial.
 
 In this paper, we study basic Dolbeault cohomology and give a new Weitzenb\"ock formula for the basic 
 Dolbeault Laplacian 
 in foliations with a transverse K\"ahler structure, and in particular we are interested in general results that 
 apply whether or not the foliation is taut.
This paper is organized as follows: In Section~\ref{Prelim}, we recall the basic properties of 
Riemannian foliations, basic cohomology,
and the basic Laplacian. 
In particular, we review the de Rham decomposition, Poincar\'e duality, the generalized Weitzenb\"ock formula and some vanishing theorems.
In Section 3,  we review the many differential operators on transversely K\"ahler foliations and 
recall some known results, providing short proofs in some cases for the reader's benefit.
If $\mathcal F$ is a minimal transversely K\"ahler foliation, then a basic holomorphic form of type $(r,0)$ is basic-harmonic (Corollary~\ref{basicholo(r,0)}). 
We prove a generalized Dolbeault decomposition (Theorem \ref{DolbeaultDecomp}). 
Also, on a transversely K\"ahler foliation, the Serre duality theorem does not necessarily hold unless $\mathcal F$ is minimal,
but we exhibit a version of Kodaira-Serre duality that
actually does hold in all cases (Theorem~\ref{Serre}, Corollary~\ref{SerreCor}, and Example~\ref{SerreFails}). 
 In Section 4, we establish some new Weitzenb\"ock formulas (Proposition~\ref{WeitzenbockKahler}) for the basic Dolbeault Laplacians
 and prove some vanishing results (Theorem~\ref{VanishingThm} and Corollary~\ref{VanishingCor}).

\section{Preliminaries}\label{Prelim}

 Let $(M,g_Q,\mathcal F)$ be a $(p+q)$-dimensional Riemannian foliation
of codimension $q$. Here, $g_Q$ is a holonomy invariant metric on the normal bundle $Q=TM/T\mathcal F$, meaning that 
$\mathcal{L}_Xg_Q=0$ for all $X\in T\mathcal F$, where $\mathcal{L}_X$ denotes the Lie derivative with respect to $X$.  
Next, let $g_M$ be a bundle-like metric on $M$ adapted to $g_Q$. This means that if $T\mathcal{F}^\perp$ is the $g_M$-orthogonal complement to 
$T\mathcal{F}$ in $TM$ and $\sigma:Q\to T\mathcal{F}^\perp$ is the canonical bundle isomorphism, then 
$g_Q=\sigma^* \left( \left. g_M\right| _{T\mathcal F^\perp}\right)$.

 Let $\pi:TM\to Q$ be the bundle projection. 
The transversal Levi-Civita
connection $\nabla$ on $Q\to M$ is given for all $s\in\Gamma(Q)$ by
\begin{equation}\label{connection}
\nabla_X s=\left\{\begin{split}& \pi([X,Y_s])\qquad \forall
X\in \Gamma(T\mathcal F)\\
 &\pi(\nabla^M_X Y_s)\qquad \forall X\in \Gamma(T\mathcal F^\perp),
\end{split}
\right.
\end{equation}
where   $Y_s \in \Gamma(T\mathcal F^\perp) $ corresponding to $s$ under the canonical isomorphism $\sigma$.  Let $R^\nabla$ and ${\rm Ric}^Q$  be the curvature tensor and
the transversal Ricci operator of $\mathcal F$, respectively.
The {\it mean curvature form} $\kappa$
of $\mathcal F$ is given by
\begin{equation}
\kappa(X)=g_Q(\sum_{i=1}^p\pi(\nabla^M_{E_i}E_i),\pi(X)) \quad \forall X\in \Gamma (TM),
\end{equation}
where $\{E_i\}_{i=1,\cdots,p}$ is a local orthonormal basis of $T\mathcal F$.
 Let $\Omega_B^r(\mathcal F)$ be
the space of all {\it basic $r$-forms}, i.e.,  $\phi\in\Omega_B^r(\mathcal F)$ if and only if
$X\lrcorner\,\phi=0$ and $\mathcal{L}_X\phi=0$ for any $X\in T\mathcal F$, where $X\lrcorner $ denotes the interior product.
$\mathcal F$ is said to be {\it minimal}  if $\kappa=0$. It is well-known(\cite{Al}) that on
a compact manifold, $\kappa_B$ is closed, i.e., $d\kappa_B=0$, where $\kappa_B$ is the basic component of $\kappa$ under the 
$L^2$-orthogonal decomposition $\Omega^1(M)=\Omega_B^1(\mathcal{F})\oplus \left(\Omega_B^1(\mathcal{F})\right)^\perp$ (\cite{Al}, \cite{PaRi}). We note that the proofs of these facts in \cite{PaRi} are still valid on a manifold 
that is not necessarily compact if the leaf closures of the Riemannian foliation are compact.
Now we recall the transversal star operator $\bar *:\Omega_B^r(\mathcal F)\to \Omega_B^{q-r}(\mathcal F)$ given by
\begin{align}
\bar * \phi = (-1)^{p(q-r)}*(\phi\wedge\chi_{\mathcal F}),
\end{align}
where $\chi_{\mathcal F}=E_1^*\wedge\ldots\wedge E_p^*$ is the characteristic form of $\mathcal F$ and $*$ is the Hodge star operator associated to $g_M$. Clearly, $\bar *^2\phi = (-1)^{r(q-r)}\phi$ for any $\phi\in\Omega_B^r\mathcal F)$. Let $\nu$ be the transversal volume form; that is, $*\nu=\chi_{\mathcal F}$. Then  the pointwise inner product $\langle\cdot,\cdot\rangle$ on $\Lambda^rQ^*$ is defined by $\langle\phi,\psi\rangle\nu = \phi\wedge\bar *\psi$. The global inner product on $L^2(\Omega_B^r(\mathcal{F}))$ is
\begin{align}
\ll\phi,\psi\gg =\int_M\langle\phi,\psi\rangle\mu_M,
\end{align}
where $\mu_M=\nu\wedge\chi_{\mathcal F}$ is the volume form with respect to $g_M$. Let 
\begin{equation}
d_B=d|_{\Omega_B^*(\mathcal F)}, \quad d_T=d_B -\epsilon(\kappa_B), \label{dBdTFormulas}
\end{equation} 
where $\epsilon(\phi)\psi=\phi\wedge\psi$ for any $\phi\in\Omega_B^1(\mathcal F)$. 
The operator $\epsilon$ acts pointwise and of course is define on $Q^*\cong (T\mathcal F^\perp)^*$.
For $v\in Q\cong T\mathcal F^\perp$, let 
\[
v\,\lrcorner=(\epsilon(v^\flat))^*
\]  
denote the interior product.
\begin{prop} (In \cite{Al},\cite{PaRi} for the compact case)  \label{deltab}
Let $(M,\mathcal{F})$ be a Riemannian foliation of codimension $q$ such that the
leaf closures of $\mathcal{F}$ are compact.
Then the formal adjoint operators $\delta_B$ and $\delta_T$ of  $d_B$ and $d_T$ with 
respect to $\ll\cdot,\cdot\gg$ on basic forms are given by 
\begin{equation*}
\delta_B\phi=(-1)^{q(r+1)+1}\bar *d_T \bar *\phi = \delta_T+\kappa_B^\sharp\lrcorner,\quad
\delta_T\phi=(-1)^{q(r+1)+1}\bar *d_B \bar *\phi, 
\end{equation*}
respectively, applied to basic $r$-forms $\phi$. 
\end{prop}
{\bf Proof.} The proofs needed from \cite{PaRi} carry over to this slightly more general case. $\Box$

Now we define two Laplacians $\Delta_B$ and $\Delta_T$ 
acting on $\Omega_B^*(\mathcal F)$ by 
\begin{equation}
\Delta_B=d_B\delta_B+\delta_B d_B, \quad \Delta_T =d_T\delta_T +\delta_T d_T,
\end{equation}
respectively. Then $\Delta_T \bar * = \bar *\Delta_B$. 
The first Laplacian $\Delta_B$ is said to be the {\it basic Laplacian}.

Now, we put $\mathcal H_B^r=\ker\Delta_B$ and $\mathcal H_T^r =\ker\Delta_T$. Then we have the following generalization of the usual De Rham-Hodge decompositon.
\begin{thm}\label{HodgeThm} (\cite{Al}, \cite{EKH}, \cite{KT3}, \cite{PaRi}) Let $(M,g_M,\mathcal F)$ be a closed oriented manifold  with a transversally oriented Riemannian foliation $\mathcal F$ and a bundle-like metric $g_M$. Then
\begin{align}
\Omega_B^r(\mathcal F) &\cong \mathcal H_B^r \oplus {\rm Im}d_B \oplus{\rm Im}\delta_B\\
&\cong \mathcal H_T^r  \oplus {\rm Im}d_T \oplus {\rm Im}\delta_T
\end{align} 
with finite dimensional $\mathcal H_B^r$ and $\mathcal H_T^r$.
\end{thm}
Since $d_B^2 = d_T^2 =0$, we define basic cohomology groups $H_B^r$ and $H_T^r$ respectively by
 $H_B^r={\ker d_B\over {\rm Im}d_B}$ and $H_T^r ={\ker d_T\over {\rm Im}d_T}$.
 Then  $\mathcal H_B^r \cong H_B^r$ and $\mathcal H_T^r\cong H_T^r$ \cite{TO2}.
  Since $\bar *\Delta_B =\Delta_T\bar*$, we have the twisted duality \cite{KT2}
\begin{align}
\mathcal H_B^r \cong \mathcal H_T^{q-r}.
\end{align}

\begin{rem}
Observe that, from the preceding results, the dimensions of $\mathcal H_B^r$ and $\mathcal H_T^r$ are smooth invariants of $\mathcal{F}$ and do not depend on the choices of bundle-like metric $g_M$ or on the transversal metric $g_Q$, even though the spaces of harmonic forms do depend on these choices.
\end{rem}

If $\mathcal F$ is a taut foliation (i.e. $\kappa=0$), then Poincar\'e duality holds. That is, $ \mathcal H_B^r \cong \mathcal H_B^{q-r}$ \cite{KT1}.
 Let $\{E_a\}$ and $\{\theta^a\}$ be the local orthonormal normal frame and its dual, respectively. 
Then we have the following  generalized Weitzenb\"ock formula.
\begin{thm} (\cite[Formula 3.9]{JU1} for the case of an isoparametric foliation)\label{WeitzThm} Let $(M,g_M,\mathcal F)$ be a 
Riemannian manifold with a foliation of codimension $q$ with compact leaf closures 
and a bundle-like metric $g_M$.
Then the generalized Weitzenb\"ock formula is given by the following:  for any basic form $\phi\in
\Omega_B^r(\mathcal F)$,
\begin{equation}
 \Delta_B \phi = \nabla_{\rm tr}^*\nabla_{\rm tr}\phi +A_{\kappa_B^\sharp}(\phi)+ F(\phi),\label{Weitz}
\end{equation}
 where
 \begin{align}
 &\nabla_{\rm tr}^*\nabla_{\rm tr}\phi =-\sum_a \nabla_{E_a}\nabla_{E_a}\phi +\nabla_{\kappa_B^\sharp}\phi,\\
 &A_Y(\phi)=\mathcal{L}_Y\phi -\nabla_Y\phi,\quad Y\in TM,\\
 &F(\phi)=\sum_{a,b}\theta^a \wedge E_b\lrcorner R^\nabla(E_b, E_a)\phi.
 \end{align}
 In particular, if $\phi$ is a basic 1-form, then $F(\phi)^\sharp={\rm Ric}^Q(\phi^\sharp)$.
\end{thm}
 It is well-known (\cite[Proposition 3.1]{JU1}) that the operator $\nabla_{\rm tr}^*\nabla_{\rm tr}$ is non-negative and
formally self-adjoint on closed manifolds, that is,
\begin{equation}
\int_M \langle\nabla_{\rm tr}^*\nabla_{\rm tr}\phi,\psi\rangle\mu_M=\int_M \langle\nabla_{\rm tr}\phi,\nabla_{\rm
tr}\psi\rangle\mu_M
\end{equation}
for all $\phi,\psi\in \Omega_B^r(\mathcal F)$, where
$\langle\nabla_{\rm tr}\phi,\nabla_{\rm tr}\psi\rangle=\sum_{a=1}^q\langle\nabla_{E_a}\phi,\nabla_{E_a}\psi\rangle$.
\begin{thm} \label{FVanishingThm}
Let $(M,g_M,\mathcal F)$ be a closed Riemannian manifold with a foliation $\mathcal F$ of codimension $q$ and a bundle-like metric $g_M$.
If the endomorphism $F$ is non-negative and positive at some point, then 
\begin{align}
\mathcal H_B^r =\{0\}.
\end{align}
\end{thm}
{\bf Proof.} Now, we modify the original bundle-like metric $g_M$ without changing the transversal metric $g_Q$ such that $\delta_B\kappa_B=0$. The existence of such metric is assured, from \cite{MO},\cite{MA}. With this new metric, the hypothesis on $F$ still holds, since this operator only depends on $g_Q$. Let $\phi$ be a basic-harmonic $r$-form. 
From (\ref{Weitz}), we have
\begin{align}
\int_M |\nabla_{\rm tr}\phi|^2 + \int_M \langle F(\phi),\phi\rangle +\int_M \langle A_{\kappa_B^\sharp}(\phi),\phi\rangle =0.
\label{Weitz=0}
\end{align}
Since $A_{\kappa_B^\sharp}(\phi) =d_B \kappa_B^\sharp\lrcorner\,\phi -\nabla_{\kappa_B^\sharp}\phi$ and $\delta_B\phi=0$, we have
\begin{align*}
\int_M \langle A_{\kappa_B^\sharp}(\phi),\phi\rangle = \int_M\langle \kappa_B^\sharp\lrcorner\,\phi,\delta_B\phi\rangle +\frac12\int_M \kappa_B^\sharp(|\phi|^2) =0.
\end{align*}
Hence from (\ref{Weitz=0}),  the proof is complete. $\Box$

Let $\mathcal R^\nabla :\Lambda^2 Q^* \to \Lambda^2 Q^*$ be the
{\it transversal curvature operator}, which is defined by
\begin{align}\label{eq1-19}
\langle\mathcal R^\nabla
(\omega_1\wedge\omega_2),\omega_3\wedge\omega_4\rangle=g_Q(R^\nabla(\omega_1^\sharp,\omega_2^\sharp)\omega_4^\sharp,\omega_3^\sharp),
\end{align}
where $\omega_i\in Q^*(i=1,\cdots,4)$.  
  Then we recall a generalization of the Gallot-Meyer theorem for foliations.
\begin{thm} \cite[Corollary 2.5]{JR} \label{RThm}Let $(M,g_M,\mathcal F)$ be as in Theorem~\ref{WeitzThm}. If $\mathcal R^\nabla\geq C$ for a positive constant $C$, then 
\begin{align}\label{eq1-22}
\langle F(\phi),\phi\rangle\geq r(q-r)C|\phi|^2
\end{align}
for any  basic
$r$-form $\phi$ $(1\leq r\leq q-1)$.
\end{thm}
From Theorem~\ref{RThm}, if the transversal curvature operator $\mathcal R^\nabla$ is positive-definite, then $F$ is also positive-definite. But the converse is not true. Hence we obtain the well-known vanishing theorem of basic-harmonic forms. 
\begin{coro} \cite[Corollary D]{MR} Let $(M,g_M,\mathcal F)$ be as in Theorem~\ref{FVanishingThm}. If the transversal curvature operator is positive-definite, then $\mathcal H_B^r =0$ for $0<r<q$.
\end{coro} 
Now, we recall the transversal divergence theorem  on a foliated manifold for later use.
\begin{thm} \cite[Theorem A]{YT}\label{TransDivThm} Let $(M,g_M,\mathcal F)$ be a closed, oriented Riemannian manifold with a transversally oriented foliation $\mathcal F$ and a bundle-like metric $g_M$ with respect to $\mathcal F$. Then
\begin{align}
\int_M {\rm div}_\nabla (\pi(X))\mu_M =\int_M g_Q(\pi(X),\kappa^\sharp)\mu_M
\end{align}
for any $X\in \Gamma (TM)$, where ${\rm div}_\nabla(\pi(X))$ is the transversal divergence of $X$ with respect to $\nabla$.
\end{thm}

\section{The basic Dolbeault cohomology on transversely K\"ahler foliations}

In this section, we generally use the same notation and cite elementary results from \cite[Section 3]{JJ} and \cite{NT}.
Let $(M,g_M,\mathcal F,J)$ be a Riemannian manifold with a transversely
K\"ahler foliation $\mathcal F$ of codimension $q=2n$.  
That is, $g_M$ is a
bundle-like metric, and $J$ is a holonomy-invariant almost complex structure 
on $Q$ such that $\nabla g_Q=0$ and $\nabla J=0$,
with $\nabla$ being the transversal Levi-Civita connection on $Q$, extended in the usual way to tensors \cite{NT}. 
In some of what follows, we will merely need the fact that the foliation is transversally Hermitian
(all of the above, merely requiring $J$ is integrable and not that $\nabla J=0$), and other times we will 
need the full power of the K\"ahler condition $\nabla J=0$.

Note that for any $X, Y\in TM$,
\begin{equation}\label{3-1}
\omega(X,Y)=g_Q(\pi(X),J\pi(Y))
\end{equation}
defines a basic K\"ahler 2-form $\omega$, which is necessarily closed.
Locally, the basic K\"ahler 2-form $\omega$ may be expressed by
\begin{align}
\omega = -\frac12\sum_{\alpha=1}^{2n}\theta^\alpha\wedge J\theta^\alpha,
\end{align}
where $\{\theta^\alpha\}_{\alpha=1,\cdots,2n}$ is a local orthonormal frame of $Q^*$.
Here, we extend $J$ to elements of $Q^*$ by setting
\begin{equation}
(J\theta)(X)=-\theta(JX)   \label{JonForms}
\end{equation}
for any $X\in Q_x,\theta\in Q_x^*$. Note that if $J^*$ denotes the dual transformation on $Q^*$, the $J$ above 
satisfies $J^*=-J$.
We prefer this convention so that $J$ commutes with the musical isomorphisms, so for example
\[
(Jv)^\flat = J(v^\flat )
\]
for all $v\in Q$.
When it is convenient, we will also refer to the bundle map $J':TM\to TM$ 
defined by $J'(v)=J(\pi(v))$ and 
abuse notation by denoting $J=J'$. Similarly, we sometimes will act on all of 
$T^*M$ using the symbol $J$.
We note that all of the above is true for transversally Hermitian foliations, but the form $\omega$ is not closed
unless the foliation is K\"ahler.

\begin{ex}
\label{KaehlerExact}
Note that in contrast to the situation of a K\"ahler form on an ordinary manifold, 
it is possible that $\omega$ is a trivial class in 
basic cohomology. 
We consider the Carri\`{e}re example from \cite%
{Ca}. Let $A$ be a matrix in $\mathrm{SL}_{2}(%
\mathbb{Z})$ of trace strictly greater than $2$. We denote respectively by $%
V_{1}$ and $V_{2}$ unit eigenvectors associated with the eigenvalues $\lambda 
$ and $\frac{1}{\lambda }$ of $A$ with $\lambda >1$ irrational. Let the
hyperbolic torus $\mathbb{T}_{A}^{3}$ be the quotient of $\mathbb{T}%
^{2}\times \mathbb{R}$ by the equivalence relation which identifies $(m,t)$
to $(A(m),t+1)$. The flow generated by the vector field $V_{2}$ is a
Riemannian foliation with bundle-like metric (letting $\left( x,s,t\right) $ denote
local coordinates in the $V_{2}$ direction, $V_{1}$ direction, and $\mathbb{R%
}$ direction, respectively)
\begin{equation*}
g=\lambda ^{-2t}dx^{2}+\lambda ^{2t}ds^{2}+dt^{2}.
\end{equation*}%
Note
that the mean curvature of the flow is $\kappa =\kappa _{b}=\log \left(
\lambda \right) dt$, since $\chi _{\mathcal{F}}=\lambda ^{-t}dx$ is the
characteristic form and $d\chi _{\mathcal{F}}=-\log \left( \lambda \right)
\lambda ^{-t}dt\wedge dx=-\kappa \wedge \chi _{\mathcal{F}}$. Then an orthonormal
frame field for this manifold is $\{X=\lambda^t\partial_x, 
Y_1=\lambda^{-t}\partial_s, 
Y_2=\partial_t\}$
corresponding to the orthonormal coframe
$\{
X^*=\chi_{\mathcal{F}}=\lambda^{-t}dx,
Y_1^*=\lambda^tds,
Y_2^*=dt
\}$. Then, letting $J$ be defined by $J(Y_1)=Y_2,J(Y_2)=-Y_1$, the Nijenhuis tensor 
\[
N_J(Y_1,Y_2) = [Y_1,Y_2]+J\left( [JY_1,Y_2]+[Y_1,JY_2] \right) - [JY_1,JY_2]
\]
clearly vanishes, so that $J$ is integrable. 

The corresponding transverse K\"ahler form is
seen to be $\omega=Y_2^*\wedge Y_1^*=\lambda^tdt\wedge ds=d(\frac 1{\log\lambda}Y_1^*)$,
an exact form in basic cohomology.
\end{ex}

Let $Q^C=Q\otimes\mathbb C$ be the complexified normal bundle and let 
\begin{align*}
Q^{1,0}=\{Z\in Q^C|JZ=iZ\},\quad Q^{0,1}=\{Z\in Q^C | JZ=-iZ\}.
\end{align*}
The element of $Q^{1,0}$ (resp. $Q^{0,1}$) is called a {\it complex normal vector field of type} $(1,0)$ (resp. (0,1)).
Then $Q^C = Q^{1,0} \oplus Q^{0,1}$ and 
\begin{align}
Q^{1,0}=\{X-iJX|\ X\in Q\},\quad Q^{0,1}=\{X+iJX |\ X\in Q\}.
\end{align}
Now, let $Q^*_C$ be the real dual bundle of $Q^C$, defined at each $x\in M$ to be the $\mathbb{C}$-linear maps 
from $Q_x^C$ to $\mathbb{C}$. Let $\Lambda_CQ^*=\Lambda Q_C^*$.
We define two sub-bundles  $Q_{1,0}$ and $Q_{0,1}$ of $\Lambda_C^1 Q^*$ by
\begin{align*}
Q_{1,0}&=\{\xi\in\Lambda^1_C Q^*|\ \xi(Z) =0,\ \forall Z\in Q^{0,1}\},\\
Q_{0,1}&=\{\xi\in \Lambda^1_C Q^* |\ \xi(Z) =0,\ \forall Z\in Q^{1,0}\}.
\end{align*}
Then $\Lambda_C^1 Q^* =Q_{1,0}\oplus Q_{0,1}$ and
\begin{align*}
Q_{1,0}&=\{\theta+i J\theta |\ \theta\in Q_C^*,
\left.\theta\right|_Q\in Q^*\}\quad \text{and}\quad \\
Q_{0,1}&=\{\theta-iJ\theta |\ \theta\in Q_C^*,
\left.\theta\right|_Q\in Q^*\},
\end{align*}
where $(J\theta)(X)=-\theta(JX)$ for any $X\in Q^C$. 
Let $\Lambda^{r,s}_CQ^*$ be the subspace of $\Lambda_C Q^*$ spanned by $\xi\wedge\eta$, where $\xi\in\Lambda^r Q_{1,0}$ and $\eta\in\Lambda^s Q_{0,1}$. The sections of $\Lambda^{r,s}_CQ^*$ are said to be  {\it forms of type $(r,s)$}. Let $\Omega_B^{r,s}(\mathcal F)$ be the set of the basic forms of type $(r,s)$.
Let $\{E_a,JE_a\}_{a=1,\cdots,n}$ be a local orthonormal basic frame(called a {\it $J$-basic frame}) of $Q$ and $\{\theta^a,J\theta^a\}_{a=1,\cdots,n}$ be their dual basic forms on $Q^*$.  Let $V_a ={1\over\sqrt 2}(E_a -iJE_a)$ and $\omega^a ={1\over\sqrt 2}(\theta^a +iJ\theta^a) $. Then 
\begin{align*}
\omega^a(V_b)=\bar\omega^a(\bar V_b)=\delta_{ab},\ \omega^a(\bar V_b)=\bar\omega^a(V_b)=0.
\end{align*}
A frame field $\{V_a\}$ is a local orthonormal basic frame field on $Q^{1,0}$, which is called a {\it normal frame field of type $(1,0)$}, and $\{\omega^a\}$ is a dual coframe  field on $Q_{1,0}$. 

In what follows, we extend the connection $\nabla$ on $Q$ in the natural way so that $\nabla_XY$ is defined for any $X\in\Gamma (TM\otimes \mathbb{C})$ and any $Y\in \Gamma (Q^C)$.
We further extend it to differential forms, requiring that $\nabla$ is a Hermitian connection, i.e., for any $V\in Q^C$ and any $\phi,\psi\in\Omega_B^{r,s}(\mathcal{F})$,
\begin{align}
V\langle\phi,\psi\rangle =\langle\nabla_V\phi,\psi\rangle +\langle\phi,\nabla_{\bar V}\psi\rangle.
\end{align}
It is an easy exercise to show that for any complex vector field $X$, $\nabla_X$ preserves the $(r,s)$ type of the form or vector field.

\begin{defn} {\rm Let $(M,g_Q,\mathcal F,J)$ be a transversely holomorphic foliation. Then  a complex normal vector field  $Z$ of type (1,0) is said to be} transversally holomorphic {\rm if 
\begin{align}
\nabla_{\bar V_a}Z=0\text{ for }1\le a \le n.
\end{align} }
\end{defn}
Let $V(\mathcal F)$  be the set of infinitesimal automorphisms $X$, i.e., $[X,Y]\in T\mathcal F$ for any $Y\in T\mathcal F$. 
The Lie derivative $\mathcal{L}_X$ can be applied to tensors with inputs in $Q$ or $Q^*$, because all such tensors can be extended by $0$ to tensors with inputs in $TM$ and $T^*M$.
We will now consider what is true when the K\"ahler condition $\nabla J=0$ holds.

\begin{prop} \label{equivProp}Let $(M,g_Q,\mathcal F,J)$ be a transversely K\"ahler foliation with compact leaf closures, with bundle-like metric $g_M$ and $X\in V(\mathcal F)$. 
Then the following are equivalent.

$(1)$  A vector field $\pi(X)$ is transversally automorphic, i.e., $\mathcal{L}_XJ=0$ $($if and only if 
$\nabla_{JY}\pi(X) = J\nabla_Y\pi(X)$ for any $Y\in T\mathcal F^\perp$ $)$.

$(2)$ A complex normal vector field  $Z=\pi(X)-iJ\pi(X)$ is transversally holomorphic, i.e., $\nabla_{\bar V_a}Z=0$ for $1\le a \le n$.
\end{prop}
Now we recall the following proposition.
\begin{prop} \cite[Theorem 5.2]{JL} and \cite[Corollary 5.10]{JJ} Let $(M,g_Q,\mathcal F,J)$  be a transversely K\"ahler foliation on a compact Riemannian manifold
with bundle-like metric $g_M$.  Then
for any $X\in V(\mathcal F)$,

$(1)$ $\pi(X)$ is transversally automorphic if and only if
\begin{align*}
&\nabla_{\rm tr}^*\nabla_{\rm tr}\pi(X) -{\rm Ric}^Q(X) + A_{X}\kappa_B^\sharp =0,\\
 &\int_M g_Q((\mathcal{L}_XJ)\kappa_B^\sharp,J\pi(X)) =0.
 \end{align*}
 
$(2)$ $\pi(X)$ is transversally Killing $($i.e., $\mathcal{L}_Xg_Q=0)$ if and only if 
\begin{align*}
\nabla_{\rm tr}^*\nabla_{\rm tr}\pi(X) - {\rm Ric}^Q(X) + A_{X}\kappa_B^\sharp =0,\quad
{\rm div}_\nabla(\pi(X))=0.
\end{align*}
\end{prop}

\begin{rem} {\rm Let $(M,g_Q,\mathcal F,J)$  be a transversely K\"ahler foliation on a compact Riemannian manifold
with bundle-like metric $g_M$. Then from \cite[Theorem 5.8, Corollaries 5.9 and 5.10]{JJ}

 (1)  If $\mathcal F$ is minimal, then every transversal Killing  field is transversally automorphic. 

 (2)  If $\mathcal F$ is non-minimal and the transversal Ricci curvature is positive, then every transversal Killing field is transversally automorphic.}
\end{rem}
Now, let $\langle\cdot,\cdot\rangle$ be a Hermitian inner product on $\Lambda_C^{r,s}(\mathcal F)$ induced by the metric, and let $\bar*:\Lambda_C^{r,s}(\mathcal F)\to \Lambda_C^{n-s,n-r}(\mathcal F)$ be the star operator defined by 
\begin{align}
\phi\wedge\bar *\bar\psi = \langle\phi,\psi\rangle \nu, 
\end{align}
where $\nu={\omega^n\over n!}$ is the transverse volume form. 

Then for any $\psi\in\Lambda_C^{r,s}(\mathcal F)$,
\begin{align}
\overline{\bar*\psi} =\bar*\bar\psi,\quad \bar*^2 \psi = (-1)^{r+s}\psi.
\end{align}

We define  the operators $\partial_B:\Omega_B^{r,s}(\mathcal F)\to \Omega_B^{r+1,s}(\mathcal F)$ and $\bar\partial_B:\Omega_B^{r,s}(\mathcal F)\to \Omega_B^{r,s+1}(\mathcal F)$  in the usual way using the differential 
$d=\partial_B+\bar \partial_B$
and the projections $\Omega_B\to\Omega_B^{r,s}$. Because of the given transverse Hermitian structure, the complex differentials satisfy 
\begin{align}
\partial_B\phi &=\sum_{a=1}^n \omega^a\wedge\nabla_{V_a}\phi,\\
\bar\partial_B\phi&=\sum_{a=1}^n\bar\omega^a\wedge\nabla_{\bar V_a}\phi, \label{partialBForm}
\end{align}
respectively. 

We now write
\begin{eqnarray}
\kappa_B &=& \kappa_B^{1,0}+\kappa_B^{0,1} \notag \\
 \kappa_B^{1,0} &=& \frac 12 (\kappa_B + i J\kappa_B ) \in\Omega_B^{1,0}(\mathcal{F})   \label{kappaComponents}\\
 \kappa_B^{0,1} &=&  \overline{\kappa_B^{1,0}}  = \frac 12 (\kappa_B - i J\kappa_B ) \in\Omega_B^{0,1}(\mathcal{F}) \notag
\end{eqnarray}

Now we utilize the bundle-like metric $g_M$ to define the operators $\partial_T:\Omega_B^{r,s}(\mathcal F)\to \Omega_B^{r+1,s}(\mathcal F)$ and $\bar\partial_T:\Omega_B^{r,s}(\mathcal F)\to\Omega_B^{r,s+1}(\mathcal F)$ respectively by
\begin{align}
\partial_T\phi =\partial_B\phi -\epsilon(\kappa_B^{1,0})\phi,\quad \bar\partial_T\phi =\bar\partial_B\phi -\epsilon(\kappa_B^{0,1})\phi,
\end{align}
Then by (\ref{dBdTFormulas}) we have
\begin{align}
d_T =\partial_T +\bar\partial_T.
\end{align}
From Proposition~\ref{deltab}, the adjoint operators $\delta_B $ and $\delta_T$ of $d_B$ and $d_T$ are given  by
\begin{align}
\delta_B = - \bar* d_T \bar*, \quad\delta_T = -\bar * d_B\bar *,\label{deltabtformulas}
\end{align}
respectively. 

  Let  $\partial_T^*$, $\bar\partial_T^*$, $\partial_B^*$ and $\bar\partial_B^*$ be the formal adjoint operators of $\partial_T,\ \bar\partial_T,\ \partial_B$ and $\bar\partial_B$, respectively, on the space of basic forms. Then 
  \begin{align}
  \delta_B =\partial_B^* +\bar\partial_B^*,\quad \delta_T =\partial_T^* +\bar\partial_T^*,
  \end{align}
  and from (\ref{deltabtformulas}),
\begin{align}
&\partial_T^*\phi = -\bar *\bar\partial_B\bar *\phi,\quad \bar\partial_T^*\phi =-\bar *\partial_B\bar *\phi, \label{starBarDeltaForm1} \\
&\partial_B^* \phi = -\bar * \bar\partial_T \bar *\phi,\quad \bar\partial_B^* \phi =-\bar *\partial_T\bar*\phi. \label{starBarDeltaForm2} 
\end{align}
These formulas were first used in the case of minimalizable (homologically oriented) transverse K\"ahler foliations in 
\cite[Section 3.4]{EK}.
In what follows,
we let
\begin{equation}
H^{1,0}=\overline{(\kappa_B^{1,0})^\sharp}=\frac 12(\kappa_B^\sharp-iJ\kappa_B^\sharp),\quad 
H^{0,1} = \overline{H^{1,0}}. \label{H10Defn}
\end{equation}
We will extend the definitions
of exterior product $\epsilon$ and interior product $\lrcorner$ linearly to complex vectors. Observe that $V\lrcorner$ is by definition the adjoint of 
$\epsilon(V^\flat)$ on real vector fields, but on complex vector fields an adjustment must be made.
If $v,w$ are real tangent vectors, the interior product satisfies
\begin{eqnarray*}
(v-iw)\lrcorner  &= &(v\lrcorner ) - i (w\lrcorner)
= (\epsilon(v^\flat)-i\epsilon(w^\flat))^*
 \\
&=& (\epsilon(v^\flat-iw^\flat))^*. 
\end{eqnarray*}
So if we denote  $X^\flat = v^\flat + iw^\flat$ for any $X = v + iw$, then  for complex vectors $X$,
\begin{equation}
(X\lrcorner )^*=\epsilon( X^\flat), \quad(\epsilon(X^\flat))^*= X\lrcorner. \label{cpxIntProd}
\end{equation}
Then for complex vectors $X$ it follows that 
\begin{equation}
\bar *\epsilon(X^\flat)\bar* =  X\lrcorner, \quad \bar *(X\lrcorner)\bar* =- \epsilon( X^\flat)
\label{starBarIntProd}
\end{equation}
Then we have the following.
\begin{prop} \label{partialBadjointFormulasProp}
(Can be deduced from \cite[Section 3]{JJ} and \cite[Theorem 3.1]{JP} in transversely K\"ahler case) If $(M,g_Q,\mathcal F,J)$  is a transversely Hermitian foliation with compact leaf closures
and bundle-like metric $g_M$, we have
\begin{align}
\partial^*_B \phi =\partial_T^*\phi +H^{1,0}\lrcorner\,\phi,\quad 
\bar\partial_B^*\phi=\bar\partial_T^*\phi +H^{0,1}\lrcorner\,\phi,
\end{align} 
where
\begin{align}
\partial_T^*\phi = -\sum_{a=1}^n V_a\lrcorner \nabla_{\bar V_a}\phi,\quad \bar\partial_T^*\phi=  -\sum_{a=1}^n \bar V_a\lrcorner \nabla_{V_a}\phi.
\end{align} 
\end{prop}
{\bf Proof.} Observe that we just apply the operators  $\delta_B^*$ and $\delta_T^*$ to 
forms of type $(r,s)$ and take the $(r-1,s)$ and $(r,s-1)$ components of the result, using
the standard formulas for the divergence and Proposition~\ref{deltab}. $\Box$

We note that similar calculations were done for the twisted Dolbeault operator in \cite[Section 3]{HV}.

We now assume that our foliation is transversally K\"ahler.
Let $L:\Omega_B^r(\mathcal F)\to \Omega_B^{r+2}(\mathcal F)$ and $\Lambda:\Omega_B^r(\mathcal F)\to \Omega_B^{r-2}(\mathcal F)$ be given by
\begin{align}
L(\phi)=\omega\wedge\phi,\quad\Lambda(\phi)=\omega\lrcorner\,\phi,
\end{align}
respectively, where $(\xi_1 \wedge\xi_2)\lrcorner =\xi_2^\sharp\lrcorner \xi_1^\sharp\lrcorner$ for any basic 1-forms $\xi_i (i=1,2)$. 
We have $\langle L\phi,\psi\rangle =\langle\phi,\Lambda\psi\rangle$ and $\Lambda=(-1)^j\bar*L\bar*$ on basic $j$-forms. For any real vector $X\in Q$, we have \cite{JJ} that 
\begin{align}
&[L,X\lrcorner ]=\epsilon(JX^b),\quad [\Lambda,\epsilon(X^b)]=- (JX)\lrcorner ,\label{LComm1} \\
&[L,\epsilon(X^b)]=[\Lambda,X\lrcorner ]=0.
\end{align}
A simple calculation shows that the formulas above extend to complex vectors $X$ as well.
Now, we extend the complex structure $J$ to $\Omega_B^r(\mathcal F)$ by 
\begin{align}
J\phi =\sum_{a=1}^{2n}J\theta^a\wedge E_a\lrcorner\,\phi.
\end{align}
Then this formula is consistent with (\ref{JonForms}), $J:\Omega_B^{r,s}(\mathcal F)\to \Omega_B^{r,s}(\mathcal F)$ is skew-Hermitian, i.e., $\langle J\phi,\psi\rangle +\langle\phi,J\psi\rangle=0$, and $J\phi = i(s-r)\phi$ for any $\phi\in\Omega_B^{r,s}(\mathcal F)$.
In addition $(J\alpha)^\sharp=J(\alpha^\sharp)$ for any one-form $\alpha$.
Note that this $J$ is much different from the operator $C$ induced from the pullback $J^*$ used often in K\"ahler geometry.

\begin{prop} \cite[Proposition 3.3]{JJ} \label{LProp}If $(M,g_Q,\mathcal F,J)$  is a transversely K\"ahler foliation 
with compact leaf closures and a bundle-like metric $g_M$, we have 
\begin{align*}
 [L,J]=[\Lambda,J]=[L,d_B]=[\Lambda,\delta_B]=0.
\end{align*}
\end{prop}
Then we have the following corollary.
\begin{coro} \cite[Proposition 3.4]{JJ}\label{LCor} With the same hypothesis, we have
\begin{align}
&[L,\partial_B]=[L,\bar\partial_B] =[\Lambda,\partial_B^*]=[\Lambda,\bar\partial_B^*]=0,\\
&[L,\partial_B^*]=-i\bar\partial_T,\ [L,\bar\partial_B^*]=i\partial_T,\ [\Lambda,\partial_B]=-i\bar\partial_T^*,\ [\Lambda,\bar\partial_B]=i\partial_T^*.
\end{align}
\end{coro}
\begin{rem}All equations in Proposition~\ref{LProp} and Corollary~\ref{LCor} hold if  we exchange the operators $(\cdot)_B$ and $(\cdot)_T$. These results were shown in the minimal foliation case, when $(\cdot)_B=(\cdot)_T$, in \cite[Lemma 3.4.4]{EK}.
\end{rem}

\bigskip
\noindent Let $\square_B$ and $\overline\square_B$
be  Laplace operators, which are defined by
\begin{align}
\square_B =\partial_B \partial_B^* + \partial_B^*\partial_B \quad{\rm and}\quad\overline\square_B =  \bar\partial_B\bar\partial_B^* + \bar\partial_B^*\bar\partial_B, \label{boxLaplaceDef}
\end{align}
respectively. Clearly, $\square_B$ and $\overline\square_B$ preserve the types of the forms.
\begin{lemma}  \label{formulaLemma}
Let $(M,g_Q,\mathcal F,J)$  be a transversely K\"ahler foliation with compact leaf closures 
and bundle-like metric $g_M$. Then
\begin{align}
\square_B &= \overline\square_B +\partial_B H^{1,0}\lrcorner + H^{1,0}\lrcorner\,\partial_B -\bar\partial_B H^{0,1}\lrcorner - H^{0,1}\lrcorner\,\bar\partial_B,   \label{LaplaceForm1}\\
\Delta_B 
&= \square_B +\overline\square_B +\partial_B H^{0,1}\lrcorner + H^{0,1}\lrcorner\,\partial_B +\bar\partial_B H^{1,0}\lrcorner + H^{1,0}\lrcorner\,\bar\partial_B. \label{LaplaceForm2}\\
0&=\bar\partial_T^*\bar\partial_B+\bar\partial_B\bar\partial_T^*-\partial_T^*\partial_B-\partial_B\partial_T^* \label{realop}
\end{align}
If $\mathcal F$ is minimal, then
\begin{align}
\Delta_B = 2\square_B = 2\overline\square_B. \label{LaplaceBoxLaplaceForm}
\end{align}
\end{lemma}
{\bf Proof.} Since $d_B^2=d_T^2=0$,  $\partial_B^2 =\bar\partial_B^2 =\partial_B\bar\partial_B +\bar\partial_B\partial_B=0$ and $\partial_T^2 =\bar\partial_T^2 =\partial_T\bar\partial_T +\bar\partial_T\partial_T=0$. 
From the equations above, Proposition~\ref{partialBadjointFormulasProp} and Corollary~\ref{LCor}, we have
\begin{align}
i(\partial_B\partial_B^* +\partial_B^*\partial_B)&=\partial_B\Lambda\bar\partial_T +\Lambda\bar\partial_T\partial_B -\partial_B\bar\partial_T\Lambda-\bar\partial_T\Lambda\partial_B\notag\\
&=[\partial_B,\Lambda]\bar\partial_B +\bar\partial_B[\partial_B,\Lambda]+\partial_B[\epsilon(\kappa_B^{0,1}),\Lambda]+[\epsilon(\kappa_B^{0,1}),\Lambda]\partial_B\notag\\
&=i(\bar\partial_T^*\bar\partial_B +\bar\partial_B\bar\partial_T^*)+\partial_B[\epsilon(\kappa_B^{0,1}),\Lambda]+[\epsilon(\kappa_B^{0,1}),\Lambda]\partial_B. \label{firstCalc}
\end{align}
Similarly, we have
\begin{align}
i(\bar\partial_B \partial_B^* + \partial_B^* \bar\partial_B) &= \bar\partial_B [\epsilon(\kappa_B^{0,1}),\Lambda] + [\epsilon(\kappa_B^{0,1}),\Lambda]\bar\partial_B,\label{Last1} \\
i(\partial_B \bar\partial_B^* + \bar\partial_B^* \partial_B) &= -\partial_B [\epsilon(\kappa_B^{1,0}),\Lambda] - [\epsilon(\kappa_B^{1,0}),\Lambda]\partial_B. \label{Last2}
\end{align}
From (\ref{LComm1}) and (\ref{H10Defn}), we have $[\epsilon(\kappa_B^{0,1}),\Lambda]=i H^{1,0}\lrcorner $ and $[\epsilon(\kappa_B^{1,0}),\Lambda]=-i H^{0,1}\lrcorner $. Hence (\ref{LaplaceForm1}) is proved from (\ref{firstCalc}) 
using Proposition~\ref{partialBadjointFormulasProp}. By direct calculation, we have
\begin{align}
\Delta_B =\square_B +\overline\square_B + (\bar\partial_B\partial_B^*+\partial_B^*\bar\partial_B) + (\partial_B\bar\partial_B^* +\bar\partial_B^*\partial_B).
\end{align}
From (\ref{Last1}) and (\ref{Last2}), the proof of (\ref{LaplaceForm2}) is complete. 
Formula (\ref{realop}) follows from adding (\ref{firstCalc}) to its conjugate and 
comparing with (\ref{LaplaceForm1}).
$\Box$

The part of the Lemma above for the minimal (homologically orientable) foliation case was first shown in \cite[Proposition 3.4.5]{EK}. 
In that case, a consequence \cite[Th\'eor\`eme 3.4.6]{EK} is that $\mathcal{H}^k_B=\bigoplus_{r+s=k}\mathcal{H}_B^{r,s}$ for all $k$, where $\mathcal H_B^{r,s}=\ker \overline\square_B$,  which is in general false.

\begin{coro} \label{CorLaplaceType} With the hypothesis of Lemma~\ref{formulaLemma}:\newline
If $\phi$ is of type $(r,0)$, then
\begin{align}
\Delta_B\phi = 2\overline\square_B\phi +
\mathcal{L}_{H^{1,0}}\phi
-H^{0,1}\lrcorner\,\bar\partial_B\phi.\label{FirstCorForm}
\end{align}
If $\phi$ is of type $(r,n)$, then
\begin{align}
\Delta_B\phi = 2\overline\square_B\phi +\partial_B \kappa_B^\sharp\lrcorner\,\phi+ \kappa_B^\sharp\lrcorner  \partial_B\phi -\bar\partial_B H^{0,1}\lrcorner\,\phi.\label{2CorForm}
\end{align}
If $\phi$ is of type $(0,s)$, then
\begin{align}
\Delta_B\phi = 2\overline\square_B\phi +(\partial_B-\bar\partial_B)H^{0,1}\lrcorner\,\phi + H^{0,1}\lrcorner (\partial_B -\bar\partial_B)\phi +H^{1,0}\lrcorner\,\partial_B\phi.
\label{neededCorForm}
\end{align}
If $\phi$ is of type $(n,s)$, then
\begin{align}
\Delta_B\phi = 2\overline\square_B\phi - i(J\kappa_B^\sharp)\lrcorner\,\bar\partial_B\phi -i\bar\partial_B(J\kappa_B^\sharp)\lrcorner\,\phi +\partial_B H^{1,0}\lrcorner\,\phi.\label{3CorForm}
\end{align}
\end{coro}
{\bf Proof.} Let $\phi$ be a type $(r,0)$. Then $H^{0,1}\lrcorner\,\phi = H^{0,1}\lrcorner\,\partial_B\phi=0$.  From (\ref{LaplaceForm1}) and (\ref{LaplaceForm2}) and Cartan's formula for the Lie derivative, the proof of (\ref{FirstCorForm}) is completed. The others are similarly proved. $\Box$

\begin{coro} Let $\phi$ be a basic form of type $(r,s)$. Then  $\phi $ is basic-harmonic   if and only if $(\square_B+\overline\square_B)\phi=0$ and
\begin{align}
 \{H^{0,1}\lrcorner\,\partial_B +\partial_B H^{0,1}\lrcorner \}\phi = \{H^{1,0}\lrcorner\,\bar\partial_B +\bar\partial_B H^{1,0}\lrcorner \}\phi=0.
\end{align}
\end{coro}
{\bf Proof.} Let $\phi\in\Omega_B^{r,s}(\mathcal F)$ be a basic form of type $(r,s)$. Clearly, $(\square_B+\overline\square_B)\phi\in\Omega_B^{r,s}(\mathcal F)$,  $\{H^{0,1}\lrcorner\,\partial_B +\partial_B H^{0,1}\lrcorner \}\phi\in\Omega_B^{r+1,s-1}(\mathcal F)$ and $\{H^{1,0}\lrcorner\,\bar\partial_B +\bar\partial_B H^{1,0}\lrcorner \}\phi\in\Omega_B^{r-1,s+1}(\mathcal F)$. Hence from (\ref{LaplaceForm2}), the proof follows. $\Box$

\begin{coro}\label{Laplacefcns}
Let $(M,g_Q,\mathcal F,J)$  be a transversely K\"ahler foliation with compact leaf closures 
and a bundle-like metric $g_M$. Then the Laplacian on complex-valued functions satisfies
\begin{align}
\square_B &= \overline\square_B -i(J\kappa_B)^\sharp,   \label{LaplaceF1fcns}\\
\Delta_B 
&= \square_B +\overline\square_B, \label{LaplaceF2fcns}\\
\bar\partial_T^*\bar\partial_B&=\partial_T^*\partial_B. \label{realopfcns}
\end{align}
The second equation holds more generally for all transversely Hermitian foliations.
\end{coro}
{\bf Proof.} 
Equation (\ref{LaplaceF1fcns}) follows from (\ref{LaplaceForm1}), using $H^{1,0}\lrcorner\,\partial_Bf=H^{1,0}\lrcorner\,(\partial_B+\bar\partial_B)f=H^{1,0}\lrcorner\, d_Bf= H^{1,0}f$ for all complex-valued functions $f$. 
Equation (\ref{LaplaceF2fcns}) comes from expanding $\Delta_B=d_B\delta_B+\delta_B d_B$ with $d_B=\partial_B+\bar\partial_B$.
Equation (\ref{realopfcns}) follows from (\ref{realop}).
 $\Box$

\begin{coro}
Let $(M,g_Q,\mathcal F,J)$  be a transversely K\"ahler foliation with compact leaf closures 
and a bundle-like metric $g_M$. Then the Laplacian on complex-valued functions satisfies $\Delta_B=2\overline\square_B$ if
and only if $\kappa_B=0$.
\end{coro}

\begin{lemma} \label{transholoLem} Let $(M,\mathcal F,J)$ be a transversely Hermitian foliation.   
Then a complex vector field $Z\in Q^{1,0}$ of type $(1,0)$ is transversally holomorphic (equivalently, $Z+\bar Z$ is transversally automorphic) if and only if 
\begin{align}
\bar\partial_B Z\lrcorner +Z\lrcorner\,\bar\partial_B=0. \label{delBarZForm}
\end{align}
\end{lemma}
{\bf Proof.} Let $Z\in Q^{1,0}$ be a transversally holomorphic vector field, that is,
$\nabla_{\bar V_a} Z =0$ ($1\le a\le n$). Then for any $\phi\in\Omega_B^{r,s}(\mathcal{F})$,
\begin{align}
\bar\partial_B Z\lrcorner\,\phi&=\sum_a\bar\omega^a\wedge (\nabla_{\bar V_a}Z)\lrcorner\,\phi +  \sum_a \bar\omega^a\wedge Z\lrcorner \nabla_{\bar V_a}\phi\\
&=-Z\lrcorner\,\bar\partial_B\phi,
\end{align}
which proves (\ref{delBarZForm}). Conversely, if $Z$ satisfies (\ref{delBarZForm}),  then
\begin{align*}
0&=\bar\partial_B Z\lrcorner\,\omega^b + Z\lrcorner\,\bar\partial_B\omega^b\\
&=\sum_a \bar\omega^a\wedge \nabla_{\bar V_a}( Z\lrcorner\,\omega^b) -\sum_a \bar\omega^a\wedge Z\lrcorner\nabla_{\bar V_a}\omega^b\\
&=\sum_a (\nabla_{\bar V_a}Z\lrcorner\,\omega^b)\bar\omega^a.
\end{align*}
Hence $\nabla_{\bar V_a}Z\lrcorner\,\omega^b=0$ for any $a$ and $b$, which means $\nabla_{\bar V_a}Z=0$ for any $a$. That is, $Z$ is transversally holomorphic.  $\Box$
\begin{defn} {\rm On a transversely holomorphic foliation, a basic form $\phi$ is said to be} basic holomorphic {\rm if $\bar\partial_B\phi=0$.}
\end{defn} 
\begin{prop} \label{kappaHoloProp}Let $(M,g_Q,\mathcal F,J)$  be a transversely K\"ahler foliation with compact
leaf closures and a bundle-like metric $g_M$, such that $\kappa_B^\sharp$ is transversally automorphic. If $\phi$ is a basic holomorphic form of type $(r,0)$, then $\square_B\phi =\Delta_B\phi$ and $(\partial_B\partial_T^* +\partial_T^*\partial_B)\phi=0$. 
\end{prop}
{\bf Proof.} 
Let $\phi$ be a basic holomorphic form of type $(r,0)$, i.e., $\bar\partial_B\phi=0$. Moreover, $\bar\partial^*_B\phi=0$ because $\bar\partial^*_B\phi$ is of type $(r,-1)$.  Hence 
 $\overline\square_B\phi=0$. Since $\kappa_B^\sharp$ is transversally automorphic, from Lemma~\ref{transholoLem}
\begin{align}
\bar\partial_BH^{1,0}\lrcorner\,\phi=0.
\end{align}
That is, $H^{1,0}\lrcorner\,\phi$ is a basic holomorphic form of type $(r-1,0)$.
 On the other hand, 
\begin{align} 
H^{0,1}\lrcorner\,\phi = H^{0,1}\lrcorner\,\partial_B\phi =0,
\end{align}
because $H^{0,1}\lrcorner\,\phi$ is of type $(r,-1)$ and $H^{0,1}\lrcorner\,\partial_B\phi$ is of type $(r+1,-1)$.
Hence from  (\ref{LaplaceForm2}), we have 
\begin{align}
 \square_B\phi=\Delta_B\phi.
\end{align}
And from (\ref{realop}), $(\partial_B\partial_T^* +\partial_T^*\partial_B)\phi=0$ is proved. $\Box$


The following example shows that in general we do not expect the mean curvature to be
transversally automorphic for any metric.
\begin{ex}
\label{KaehlerExactMeanCurv} We continue from Example \ref{KaehlerExact}.
From previous calculations, 
\begin{eqnarray*}
\kappa  &=&\kappa _{B}=\log \left( \lambda \right) dt=\log \left( \lambda
\right) Y_{2}^{\ast } \\
Y_{1}^{\ast } &=&\lambda ^{t}ds,~Z^{\ast }=\frac{1}{2}\left( Y_{1}^{\ast
}+iY_{2}^{\ast }\right) ,
\end{eqnarray*}%
so%
\begin{eqnarray*}
\kappa _{B} &=&-i\left( \log \lambda \right) Z^{\ast }+i\left( \log \lambda
\right) \bar{Z}^{\ast } \\
&=&-i\left( \log \lambda \right) \frac{1}{2}\left( \lambda ^{t}ds+idt\right)
+i\frac{1}{2}\left( \log \lambda \right) \bar{Z}^{\ast }.
\end{eqnarray*}%
Then%
\begin{eqnarray*}
\kappa _{B}^{1,0} &=&-i (\log\lambda) Z^* =
-\frac{i}{2}\left( \log \lambda \right) \left( \lambda
^{t}ds+idt\right)  \\
\bar{\partial}_{B}\kappa _{B}^{1,0} &=&d\kappa _{B}^{1,0}=\frac{i}{2}\left(
\log \lambda \right) ^{2}\lambda ^{t}ds\wedge dt \\
&=&\frac{i}{2}\left( \log \lambda \right) ^{2}Y_{1}^{\ast }\wedge
Y_{2}^{\ast } \\
&=&\left( \log \lambda \right) ^{2}\bar{Z}^{\ast }\wedge Z^{\ast }.
\end{eqnarray*}%
It is impossible to change the metric so that this is zero.
The reason is that from \cite{Al} the mean curvature $\kappa_B'$ for any
other compatible bundle-like metric would satisfy $\kappa_B'=\kappa_B+df$
for some real basic function $f$, which would imply that $(\kappa_B^{1,0})'=
\kappa_B^{1,0}+\partial_Bf$, and $\partial_Bf=Z(f)Z^*$. Since $f$ is a periodic function of 
$t$ alone, this is $\partial_Bf=-\frac i2 (\partial_tf )\,Z^*$. Then in that case
\begin{eqnarray*}
\bar\partial_B(\kappa_B^{1,0})' &=&d(\kappa_B^{1,0}-\frac i2 (\partial_tf )\,Z^*)\\
&=&d\left(-\frac{i}{2}\left( \log \lambda +\frac 12\partial_tf\right) \left( \lambda
^{t}ds+idt\right) \right)\\
&=& 
\left(\frac i2 (\log\lambda)^2 +\frac i4 (\partial_t^2 f+(\log\lambda )\partial_t f)\right)
\lambda^t ds\wedge dt \\
&=&\left((\log\lambda)^2 +\frac 12 (\partial_t^2 f+(\log\lambda )\partial_t f)\right)
\bar Z^*\wedge Z^*
\end{eqnarray*}
Since the term in parentheses is never zero for any periodic function $f$, 
we conclude that $\bar\partial \kappa_B^{1,0}$ is a nonzero multiple 
of $\bar Z^*\wedge Z^*$ for any compatible bundle-like metric.
Further, we see that 
$\bar\partial_B=\bar Z^*\wedge\nabla_{\bar Z}$, so that 
\[
\nabla_{\bar Z} \kappa_B^{1,0}=(\log\lambda)^2 Z^*,
\]
and likewise this would be nonzero for any compatible bundle-like metric. Hence $\kappa_B^\sharp$ is 
not automorphic in general for nontaut transversely K\"ahler foliations.
\end{ex}

\begin{prop}\label{(r,0)Prop}
Let $(M,g_Q,\mathcal F,J)$  be a transversely K\"ahler foliation with compact leaf closures 
and a bundle-like metric $g_M$. 
Let $\phi$ be a basic holomorphic form of type $(r,0)$. Then
\[
\Delta_B\phi = \mathcal{L}_{H^{1,0}}\phi.
\]
\end{prop}
{\bf Proof.} Let $\phi$ be a basic holomorphic form of type $(r,0)$. Trivially, $\overline\square_B\phi=0$. Hence from  (\ref{FirstCorForm}), the proof follows.
$\Box$

\begin{coro} (First proved in \cite[Th\'eor\`eme 3.4.6]{EK}) \label{basicholo(r,0)}
Let $(M,g_Q,\mathcal F,J)$  be a minimal transversely K\"ahler foliation with compact leaf closures 
on a Riemannian manifold
with bundle-like metric $g_M$. Then a basic form of type $(r,0)$ is basic-harmonic
if and only if it is basic holomorphic.\end{coro}
{\bf Proof.} Let $\phi$ be a basic form of type $(r,0)$. Since automatically 
$\bar\partial_B^*\phi=0$, $\bar\partial_B\phi = 0 $ if and only if
$\overline\square_B\phi=0$ , if and only if $\Delta_B\phi=0$ by Proposition~\ref{(r,0)Prop}.
$\Box$


\bigskip
\noindent  Note that $\bar\partial_B^2=0$. So there exists a complex
\begin{align*}
0 \overset {\bar\partial_B}\longrightarrow \cdots \overset{\bar\partial_B}\longrightarrow\Omega_B^{r,s-1}(\mathcal F)\overset{\bar\partial_B}\longrightarrow\Omega_B^{r,s}(\mathcal F)\overset{\bar\partial_B}\longrightarrow\Omega_B^{r,s+1}(\mathcal F)\overset{\bar\partial_B}\longrightarrow\cdots \overset{\bar\partial_B}\longrightarrow 0.
\end{align*} 
Hence we can define the {\it basic Dolbeault cohomology group} on a transversely holomorphic foliation by
\begin{align}
H_B^{r,s}={{\ker\bar\partial_B}\over {\rm Im} \bar\partial_B}.
\end{align}
Then we have the generalization of the Dolbeault decomposition.
\begin{thm} (Proved in \cite[Th\'eor\`eme 3.3.3]{EK})
\label{DolbeaultDecomp}Let $(M,g_Q,\mathcal F,J)$  be a transversely Hermitian foliation on a closed Riemannian manifold
with bundle-like metric $g_M$. Then
\begin{align}
\Omega_B^{r,s}(\mathcal F)\cong \mathcal H_B^{r,s} \oplus {\rm Im}\bar\partial_B \oplus {\rm Im}\bar\partial_B^*,
\end{align}
where $\mathcal H_B^{r,s}=\ker\overline\square_B$ is finite dimensional. Moreover,  $\mathcal H_B^{r,s} \cong H_B^{r,s}$.
\end{thm}
{\bf Proof.} The proof is similar to the one in Theorem~\ref{HodgeThm}. See \cite{KT3} precisely. $\Box$

\bigskip
\noindent Let $\square_T  =\partial_T \partial_T^* +\partial_T^* \partial_T$ and  $\overline\square_T  =\bar\partial_T \bar\partial_T^* +\bar\partial_T^* \bar\partial_T$. Then from (\ref{starBarDeltaForm1}) and (\ref{starBarDeltaForm2}), which only 
require the transverse Hermitian structure,
\begin{align}
\bar*\overline\square_B =\square_T \bar *,\quad \bar*\square_B =\overline\square_T \bar * .\label{starBarLaplaceForm}
\end{align}
 
\begin{thm} \label{Serre}Let $(M,g_Q,\mathcal F,J)$  be a transversely Hermitian foliation on a closed Riemannian manifold
with bundle-like metric $g_M$. Then
\begin{align}
\mathcal H_B^{r,s}\cong \mathcal H_T^{n-r,n-s},\label{SerreForm}
\end{align}
where $\mathcal H_T^{r,s} =\ker \overline\square_T$.
\end{thm}
{\bf Proof.} We define the operator $\sharp:\Omega_B^{r,s}(\mathcal F)\to \Omega_B^{n-r,n-s}(\mathcal F)$ by
\begin{align*}
\sharp\phi :=\bar *\bar\phi.
\end{align*}
From (\ref{starBarLaplaceForm}), we have $\sharp \overline\square_B\phi=\overline\square_T (\sharp\phi)$. Hence if $\overline\square_B\phi =0$, then $\overline\square_T(\sharp\phi)=0$. This means that $\sharp\phi\in \ker\overline\square_T$. Since $\bar *$ is an isomorphism, it is clear that $\sharp$ is also an isomorphism. Hence  (\ref{SerreForm}) is proved. $\Box$

Let $H_T^{r,s}$ denote the basic cohomology corresponding to the differential $\partial_T$ restricted to basic forms of type $(r,s)$. This is a type of
basic Lichnerowicz cohomology. The interested reader may consult \cite[Section 3]{Ba} and \cite[Section 3, called ``adapted cohomology'' here]{Va} for information about ordinary Lichnerowicz cohomology and \cite{Had} for the basic case. The following corollary follows directly from the appropriate analogues of the Hodge theorems Theorem~\ref{HodgeThm} and 
Theorem~\ref{DolbeaultDecomp}.
A similar theorem was proved in \cite{KT2} for ordinary basic cohomology. 

\begin{coro} \label{SerreCor}
Let $(M,g_Q,\mathcal F,J)$  be a transversely Hermitian foliation on a closed Riemannian manifold
with bundle-like metric $g_M$. Then
\[
H_B^{r,s}\cong H_T^{n-r,n-s}.
\]
\end{coro}

\bigskip
\noindent If $\mathcal F$ is minimal, then $\overline\square_B =\overline \square_T$. Hence we have the following corollary.
\begin{coro} (See reference below.) If $(M,g_Q,\mathcal F,J)$  is a minimal transversely Hermitian foliation on a closed Riemannian manifold
with bundle-like metric $g_M$, then Kodaira-Serre duality   holds:
\begin{align}
\mathcal H_B^{r,s}\cong \mathcal H_B^{n-r,n-s}. \label{KSerreForm}
\end{align}
\end{coro}

\begin{rem}
The Theorem above was first proved in \cite[Theorem 3.3.4]{EK}. In Theorem 3.4.6 of the same paper, it was also shown that if in addition the foliation is minimal and K\"ahler, $\mathcal H_B^{r,s}\cong\mathcal H_B^{s,r}$, $\mathcal H_B^k =\bigoplus_{r+s=k}\mathcal H_B^{r,s}$ and that 
$\mathcal H_B^{r,r}\ne \{ 0\}$ for all $r$.
\end{rem}

\begin{ex}\label{SerreFails}
In general, Kodaira-Serre duality (\ref{KSerreForm}) fails for a transversely K\"ahler foliation on a closed Riemannian manifold with bundle-like metric, even though the standard version holds on any compact complex manifold. We return to the transversely K\"ahler foliation of Example~\ref{KaehlerExact}, with the same notation. One can easily verify that
\begin{eqnarray*}
H^{0,0}_B=\ker\overline\partial_B^{0,0}&=&\func{span}\{ 1\} \\
H^{1,0}_B=\ker\overline\partial_B^{1,0}&=&\{ 0\}  \\
H^{0,1}_B=\frac{\ker\overline\partial_B^{0,1}}{\func{im}\overline\partial_B^{0,0}}&=&\func{span}\{ Y^*_1-iY^*_2\}  \\
H^{1,1}_B=\frac{\Omega_B^{1,1}}{\func{im}\overline\partial_B^{1,0}}&=&\{ 0\}  ,\\
\end{eqnarray*}
where the last equality is true because one can show that every element of $\Omega^{1,1}_B$ is $\overline\partial_B$-exact.
Then observe that the ordinary basic cohomology Betti numbers for this foliation are $h_B^0=h_B^1=1$, $h_B^2=0$, we see that
the basic Dolbeault Betti numbers satisfy 
\[
h_B^{0,0}=h_B^{0,1}=1,\quad h_B^{1,0}=h_B^{1,1}=0.
\]
So even though it is true that
\[
h_B^j=\sum_{p+q=j} h_B^{p,q},
\]
and the foliation is transversely K\"ahler, we also have (with $n=1$)
\[
h_B^{p,q}\ne h_B^{q,p},\quad h_B^{r,s}\ne h_B^{n-r,n-s}.
\]
The exactness of the basic K\"ahler form causes the Kodaira-Serre argument, the Lefschetz theorem, the Hodge diamond ideas to fail. Thus, for a nontaut, transversely K\"ahler foliation, it is not necessarily true that the odd basic Betti numbers are even, and the basic Dolbeault numbers do not have the same kinds of symmetries as Dolbeault cohomology on K\"ahler manifolds.
Also, this example shows that the even degree basic cohomology groups are 
not always nonzero, as is the case for ordinary cohomology for symplectic manifolds (and thus all K\"ahler manifolds).
However, these analogous results really are true in the minimal case, and in fact there is an analogue of the 
Lefschetz decomposition for basic cohomology (see \cite[Section 3.4.7]{EK}).
\end{ex}

\section{The Weitzenb\"ock formula}
Let $(M,g_Q,\mathcal F,J)$  be a transversely K\"ahler foliation 
with compact leaf closures on a Riemannian manifold
with bundle-like metric $g_M$.  We define two operators
\begin{align}
\nabla_T^*\nabla_T\phi&=-\sum_a \nabla_{V_a}\nabla_{\bar V_a}\phi + \nabla_{H^{0,1}}\phi,\\
\bar\nabla_T^*\bar\nabla_T\phi&=-\sum_a \nabla_{\bar V_a}\nabla_{V_a}\phi + \nabla_{H^{1,0}}\phi.
\end{align}
Then by direct calculation, we have
\begin{align}
\nabla_T^*\nabla_T\phi=\bar\nabla_T^*\bar\nabla_T\phi +\nabla_{H^{0,1} -H^{1,0}}\phi +\sum_a R^Q(\bar V_a,V_a)\phi
\label{DeltaTSquareForm}
\end{align}
for any basic form $\phi$. And we have the following proposition.
\begin{prop} \cite[Proposition 3.2]{JP} Let $(M,g_Q,\mathcal F,J)$  be a transversely K\"ahler foliation on a compact Riemannian manifold
with bundle-like metric $g_M$. Then the operators $ \nabla_T^*\nabla_T$ and $\bar\nabla_T^*\bar\nabla_T$ are formally self-adjoint  and positive-definite.
\end{prop}
{\bf Proof.} Fix $x\in M$ and choose a normal frame field  $\{V_a\}$ of type $(1,0)$ such that $(\nabla V_a)_x =0$.  Then for any $\phi,\psi\in\Omega_B^{r,s}(\mathcal F)$,
\begin{align*}
\ll \nabla_T^*\nabla_T\phi,\psi\gg &= \int_M \langle -\sum_a \nabla_{V_a}\nabla_{\bar V_a}\phi + \nabla_{H^{0,1}}\phi,\psi\rangle\\
&=-\int_M \sum_a V_a\langle \nabla_{\bar V_a}\phi,\psi\rangle +\int_M\sum_a \langle\nabla_{\bar V_a}\phi,\nabla_{\bar V_a}\psi\rangle +\int_M \langle \nabla_{H^{0,1}}\phi,\psi\rangle\\
&=-\int_M\langle\phi,\sum_a \nabla_{ V_a}\nabla_{\bar V_a}\psi\rangle +\int_M\langle \nabla_{H^{0,1}}\phi,\psi\rangle\\
&\quad -\int_M\sum_a V_a\langle\nabla_{\bar V_a}\phi,\psi\rangle +\int_M\sum_a \bar V_a\langle\phi,\nabla_{\bar V_a}\psi\rangle.
\end{align*}
If we take complex vector fields $Y,Z\in Q^C$ satisfying
\begin{align*}
g_Q(Y,X)= \langle\nabla_{X+iJX}\phi,\psi\rangle,\quad g_Q(Z,X)=\langle\phi,\nabla_{X+iJX}\psi\rangle
\end{align*}
for any vector field $X\in\Gamma Q$, then
 by the transversal divergence theorem (Theorem~\ref{TransDivThm}), 
\begin{align}
\int_M \sum_a V_a\langle \nabla_{\bar V_a}\phi,\psi\rangle &=\frac12\int_M {\rm div}_\nabla Y =\int_M\langle\nabla_{H^{0,1}}\phi,\psi\rangle,\\
\int_M \sum_a \bar V_a\langle \phi,\nabla_{\bar V_a}\psi\rangle &=\frac12\int_M {\rm div}_\nabla Z =\int_M\langle\phi,\nabla_{H^{0,1}}\psi\rangle.
\end{align}
From the equations above,
\begin{align*}
\ll \nabla_T^*\nabla_T\phi,\psi\gg &=\int_M \langle\nabla_{\rm T}\phi,\nabla_{\rm T}\psi\rangle\\
&=\int_M\langle\phi,-\sum_a \nabla_{V_a}\nabla_{\bar V_a}\psi\rangle 
 +\int_M\langle\phi,\nabla_{H^{0,1}}\psi\rangle     \\
  &=\ll\phi,\nabla_T^*\nabla_T\psi\gg,
\end{align*}
where $\langle\nabla_{\rm T}\phi,\nabla_{\rm T}\psi\rangle =\sum_a \langle\nabla_{\bar V_a}\phi,\nabla_{\bar V_a}\psi\rangle$. Hence the proof is complete. Others are similarly proved.  $\Box$

\begin{prop} \label{WeitzenbockKahler}
If $(M,g_Q,\mathcal F,J)$  is a transversely K\"ahler foliation
with compact leaf closures on a Riemannian manifold
with bundle-like metric $g_M$, we have for all $\phi\in\Omega_B^{r,s}(\mathcal{F})$,
\begin{align}
\overline\square_B\phi&= \nabla_T^*\nabla_T\phi + \sum_{a,b}\bar\omega^a\wedge \bar V_b\lrcorner  R^Q(V_b,\bar V_a)\phi+\sum_a \bar\omega^a\wedge (\nabla_{\bar V_a}H^{0,1})\lrcorner\,\phi,  \label{boxBar}\\
\square_B\phi&= \bar\nabla_T^*\bar\nabla_T\phi  + \sum_{a,b} \omega^a\wedge V_b\lrcorner R^Q(\bar V_b ,V_a)\phi+\sum_a\omega^a\wedge (\nabla_{V_a} H^{1,0})\lrcorner\,\phi. \label{boxAfterBoxBar}
\end{align} 
\end{prop}

\noindent
{\bf Proof.} Fix $x\in M$. 
Observe that since $(M,\mathcal{F})$ is transversely K\"ahler, a local $J$-basic frame 
$\{ E_a, JE_a\}_{a=1,\ldots,n}$ can be chosen to consist of basic fields, so that 
the corresponding $V_a$ and $\omega_a$ are also basic fields of type $(1,0)$.
Furthermore, at $x\in M$, we may make the choices so that
$(\nabla V_a)_x=0$. Then at $x$, by (\ref{partialBForm}) and Proposition~\ref{partialBadjointFormulasProp} we have
for all $\phi\in\Omega_B^{r,s}(\mathcal{F})$,
\begin{align*}
\bar\partial_B \bar\partial_B^*\phi=&-\sum_{a,b} \bar\omega^a \wedge\bar V_b\,\lrcorner\, \nabla_{\bar V_a}\nabla_{V_b}\phi +\sum_a\bar\omega^a\wedge (\nabla_{\bar V_a}H^{0,1})\,\lrcorner\, \phi \\
&+\sum_a\bar\omega^a\wedge H^{0,1}\,\lrcorner\,\nabla_{\bar V_a}\phi.
\end{align*} 
Similarly, we have
\begin{align*}
\bar\partial_B^*\bar\partial_B\phi=&-\sum_a\nabla_{V_a}\nabla_{\bar V_a}\phi 
+\nabla_{H^{0,1}}\phi 
+\sum_{a,b}\bar\omega^a\wedge \bar V_b\,\lrcorner\, \nabla_{V_b}\nabla_{\bar V_a}\phi\\
&-\sum_b \bar\omega^a\wedge H^{0,1}\,\lrcorner\, \nabla_{\bar V_a}\phi.
\end{align*}
By summing the two equation above, we have (\ref{boxBar}). The proof of (\ref{boxAfterBoxBar}) is similarly proved.  $\Box$

The analogous results of the proposition above for the case of the twisted basic 
Dolbeault Laplacians was shown in \cite[Theorem 3.1]{HV}.

\begin{prop} If $(M,g_Q,\mathcal F,J)$  is a transversely K\"ahler foliation
with compact leaf closures on a Riemannian manifold
with bundle-like metric $g_M$, the following hold:

$(1)$ If $\phi$ is a basic form of type $(r,0)$, then
\begin{align}
\overline\square_B\phi =\nabla_T^*\nabla_T\phi.\label{boxBarT}
\end{align}

$(2)$ If $\phi$ is a basic form of type $(0,s)$, then
\begin{align}
\square_B\phi =\bar\nabla_T^*\bar\nabla_T\phi. \label{conjBoxBarT}
\end{align}

$(3)$ If $\phi$ is a basic form of type $(r,n)$, then
\begin{align}
\overline\square_B\phi&=\nabla_T^*\nabla_T\phi +\sum_a R^Q(V_a,\bar V_a)\phi
 +\operatorname{div}_\nabla (H^{0,1}) \phi \label{boxBarrn1}   \\
 &= \bar\nabla_T^*\bar\nabla_T \phi -i \nabla_{J\kappa_B^\sharp}\phi +\operatorname{div}_\nabla (H^{0,1}) \phi.
  \label{boxBarrn2}  
\end{align}

$(4)$ If $\phi$ is a basic form of type $(n,s)$, then
\begin{align}
\square_B\phi &=\nabla_T^*\nabla_T \phi +i \nabla_{J\kappa_B^\sharp}\phi +\operatorname{div}_\nabla (H^{1,0})\phi.
\label{boxBarrn2Conj}
\end{align}
\end{prop}
{\bf Proof.} (1) Let $\phi$ be a basic form of type $(r,0)$. Then $\bar V_b\lrcorner R^Q(V_b,\bar V_a)\phi$ and $(\nabla_{\bar V_a}H^{0,1})\lrcorner\,\phi$ are of type $(r,-1)$, which are zero. Hence from (\ref{boxBar}), Eq.(\ref{boxBarT}) is proved. 

 (2) The proof of (\ref{conjBoxBarT}) follows from the conjugation of (\ref{boxBarT}).

(3)  Now, let $\phi$ be of type $(r,n)$. Then $\bar\omega^a\wedge \bar V_b\lrcorner\,\phi =\delta_b^a \phi$, 
so that
$\operatorname{div}_\nabla (H^{0,1}) \phi = \sum_a\bar\omega^a\wedge(\nabla_{\bar V_a}H^{0,1})\lrcorner\phi$.
Since $R^Q(V_b,\bar V_a)\phi$ is also of type $(r,n)$, we have
\begin{align}
\sum_{a,b}\omega^a\wedge V_b\lrcorner R^Q(V_b,\bar V_a)\phi = \sum_aR^Q(V_a,\bar V_a)\phi. \label{RQForm}
\end{align}
From (\ref{boxBar}) and (\ref{RQForm}), Eq. (\ref{boxBarrn1}) is proved. Statement (\ref{boxBarrn2}) follows directly from (\ref{DeltaTSquareForm}).

(4) The proof of (\ref{boxBarrn2Conj}) follows from the conjugation of (\ref{boxBarrn2}). $\Box$

\begin{prop} \label{LaplaceBDeltaTProp}
Let $(M,g_Q,\mathcal F,J)$  be a transversely K\"ahler foliation 
with compact leaf closures on a Riemannian manifold
with bundle-like metric $g_M$. Then for all $\phi\in\Omega_B^{r,s}(\mathcal F)$,
\begin{align*}
\frac12\Delta_B|\phi|^2 =& \langle \nabla_T^*\nabla_T\phi,\phi\rangle + \langle \phi,\bar\nabla_T^*\bar\nabla_T\phi\rangle -\sum_{a=1}^n\{|\nabla_{\bar V_a}\phi|^2 + |\nabla_{V_a}\phi|^2\}\\
&+\frac 12 (H^{1,0}-H^{0,1})|\phi|^2,
\end{align*}
where $|\phi|^2 =\langle \phi,\phi\rangle$.
\end{prop}
{\bf Proof.} Let $\phi\in\Omega_B^{r,s}(\mathcal F)$. Since the connection $\nabla$ is  Hermitian,   from (\ref{boxBarT}) 
\begin{align*}
\overline\square_B |\phi|^2 &= -\sum_a \nabla_{V_a}\nabla_{\bar V_a} |\phi|^2 + \nabla_{H^{0,1}}|\phi|^2\\
&= \langle \nabla_T^*\nabla_T\phi,\phi\rangle + \langle \phi,\bar\nabla_T^*\bar\nabla_T\phi\rangle -\sum_{a=1}^n\{|\nabla_{\bar V_a}\phi|^2 + |\nabla_{V_a}\phi|^2\}.
\end{align*}
From Corollary~\ref{Laplacefcns} (\ref{LaplaceF1fcns}) and (\ref{LaplaceF2fcns}), we have
\begin{align*}
\Delta_B |\phi|^2 = 2\overline\square_B |\phi|^2 -i J\kappa_B^\sharp(|\phi|^2).
\end{align*}
Since $iJ\kappa_B^\sharp=H^{0,1}-H^{1,0}$, the proof follows. $\Box$

From Proposition~\ref{LaplaceBDeltaTProp}, we have the following corollary.
\begin{coro} 
\label{Laplacer0Cor}
If $\phi$ is a basic form of type $(r,0)$, then
\begin{align*}
-\frac12\Delta_B |\phi|^2 =&-\langle\overline\square_B\phi,\phi\rangle -\langle\phi,\overline\square_B\phi\rangle+\sum_a\{|\nabla_{\bar V_a}\phi|^2 +|\nabla_{V_a}\phi|^2\}\\
& +\sum_{a}\langle \phi,R^Q(\bar V_a,V_a)\phi\rangle -\frac12\langle \nabla_{H^{1,0}-H^{0,1}}\phi,\phi\rangle-\frac12\langle\phi,\nabla_{H^{1,0}-H^{0,1}}\phi\rangle.
\end{align*}
\end{coro}
{\bf Proof.} The proof follows from (\ref{DeltaTSquareForm}) and (\ref{boxBarT}). $\Box$


\begin{coro} \label{Laplacer0Cor2}
If $\phi$ is a basic-harmonic form of type $(r,0)$, then
\begin{align*}
-\frac12\Delta_B|\phi|^2=&\sum_a\{|\nabla_{\bar V_a}\phi|^2 +|\nabla_{V_a}\phi|^2\}
 +\langle \phi,\sum_aR^Q(\bar V_a,V_a)\phi\rangle\\
 &+\frac12\sum_a\{\langle\omega^a\wedge (\nabla_{V_a}H^{1,0})\lrcorner\,\phi,\phi\rangle +\langle\phi,\omega^a\wedge (\nabla_{V_a}H^{1,0})\lrcorner\,\phi\rangle\}.
 \end{align*}
\end{coro}
{\bf Proof.}
Let $\phi$ be a basic-harmonic form of type $(r,0)$. 
From (\ref{FirstCorForm}), we have
\begin{align}
&2\overline\square_B\phi = -\partial_B H^{1,0}\lrcorner\,\phi - H^{1,0}\lrcorner\,\partial_B\phi + H^{0,1}\lrcorner\,\bar\partial_B\phi,\\
&\bar\partial_B H^{1,0}\lrcorner\,\phi + H^{1,0}\lrcorner\,\bar\partial_B\phi =0.
\end{align}
Hence from Corollary~\ref{Laplacer0Cor}, we have
\begin{align*}
-\frac12\Delta_B|\phi|^2&=\frac12\langle\partial_B H^{1,0}\lrcorner\,\phi+H^{1,0}\lrcorner\,\partial_B\phi-H^{0,1}\lrcorner\,\bar\partial_B\phi,\phi\rangle\\
&+\frac12\langle\phi,\partial_B H^{1,0}\lrcorner\,\phi+ H^{1,0}\lrcorner\,\partial_B\phi-H^{0,1}\lrcorner\,\bar\partial_B\phi\rangle\\
&-\frac12\langle\nabla_{H^{1,0}-H^{0,1}}\phi,\phi\rangle -\frac12\langle\phi,\nabla_{H^{1,0}-H^{0,1}} \phi\rangle\\
&+\sum_a\{|\nabla_{\bar V_a}\phi|^2 +|\nabla_{V_a}\phi|^2\}
 +\langle \phi,\sum_aR^Q(\bar V_a,V_a)\phi\rangle.
\end{align*}
On the other hand, since $\phi$ is of type $(r,0)$, we have
\begin{align}
&H^{0,1}\lrcorner\,\bar\partial_B\phi =\nabla_{H^{0,1}}\phi,\\
&\partial_B H^{1,0}\lrcorner\,\phi + H^{1,0}\lrcorner\,\partial_B\phi = \nabla_{H^{1,0}}\phi +\sum_a \omega^a\wedge  (\nabla_{V_a}H^{1,0})\lrcorner\,\phi.
\end{align}
Hence from these two equations, we have
\begin{align*}
-\frac12\Delta_B|\phi|^2=&\sum_a\{|\nabla_{\bar V_a}\phi|^2 +|\nabla_{V_a}\phi|^2\}
 +\langle \phi,\sum_aR^Q(\bar V_a,V_a)\phi\rangle\\
 &+\frac12\sum_a\{\langle\omega^a\wedge (\nabla_{V_a}H^{1,0})\lrcorner\,\phi,\phi\rangle +\langle\phi,\omega^a\wedge (\nabla_{V_a}H^{1,0})\lrcorner\,\phi\rangle\}.
 \end{align*}
Hence the proof is complete. $\Box$

\begin{rem}\label{10Remark} {\rm Let $\xi$ and $\eta$ be basic forms of type $(1,0)$, i.e., $\xi,\eta\in\Omega_B^{1,0}(\mathcal F)$. Then
\begin{align}
\langle R^Q(\bar V_a,V_a)\xi,\eta\rangle = \langle\xi,R^Q(\bar V_a,V_a)\eta\rangle.
\end{align}
That is, $\sum_a R^Q(\bar V_a,V_a)$ is Hermitian symmetric, and so it is diagonalized by $\{\omega^a\}$. Let $\lambda_a$ be a real eigenvalue of $\sum_a R^Q(\bar V_a, V_a)$ corresponding to $\omega^a$, i.e., $(\sum_a R^Q(\bar V_a, V_a))\omega^b =\lambda_b\omega^b$. 
Then by using the first Bianchi identity, we have  
\begin{align*}
\lambda_a &=\langle \sum_b R^Q(\bar V_b,V_b)\omega^a,\omega^a\rangle\\
& =\sum_{b=1}^n g_Q (R^Q(\bar V_b,V_b)V_a,\bar V_a)\\
&=Ric^Q(E_a,E_a).
\end{align*}
}
\end{rem}
From Remark~\ref{10Remark}, we have the following corollary.
\begin{coro} \label{AnotherLaplaceCor}
Let $\phi$ be a basic-harmonic form of type $(r,0)$. Then
\begin{align*}
-\frac12\Delta_B|\phi|^2=&\sum_a\{|\nabla_{\bar V_a}\phi|^2 +|\nabla_{V_a}\phi|^2\}
 +\sum_{i=1}^r Ric^Q(E_{a_i},E_{a_i})|\phi|^2\\
 &+\frac12\sum_a\{\langle\omega^a\wedge (\nabla_{V_a}H^{1,0})\lrcorner\,\phi,\phi\rangle +\langle\phi,\omega^a\wedge (\nabla_{V_a}H^{1,0})\lrcorner\,\phi\rangle\}.
 \end{align*}
\end{coro}
{\bf Proof.} From Remark~\ref{10Remark}, we have that for any $\phi\in\Omega_B^{r,0}(\mathcal F)$, 
\begin{align*}
\langle\phi,\sum_a R^Q(\bar V_a, V_a)\phi\rangle &=(\sum_{i=1}^r \lambda_{a_i})|\phi|^2 =\sum_{i=1}^r Ric^Q(E_{a_i},E_{a_i})|\phi|^2.
\end{align*}
From Corollary~\ref{Laplacer0Cor2}, the proof is complete. $\Box$

\begin{thm} \label{VanishingThm}
Let $\mathcal F$ be a transversely K\"ahler foliation on a closed manifold with a bundle-like metric. If the transversal Ricci curvature is nonnegative and positive at some point, then for $r>0$ there are no non-zero basic-harmonic forms of type $(r,0)$; i.e., $\mathcal H_B^{r,0}=\{0\}$.
\end{thm} 
{\bf Proof.} First, we have that for any basic form $\phi$,
\begin{align}
\sum_{a=1}^n\langle\omega^a\wedge (\nabla_{V_a}H^{1,0})\lrcorner\,\phi,\phi\rangle &=\sum_{a,b=1}^n 
g_Q(\nabla_{V_a}H^{1,0},\bar V_b)\langle V_b\lrcorner\,\phi,V_a\lrcorner\,\phi\rangle\notag\\
&=\sum_{a=1}^n g_Q(\nabla_{V_a}H^{1,0},\bar V_a) |V_a\lrcorner\,\phi|^2.
\end{align}
Similarly, we have
\begin{align}
\sum_{a=1}^n\langle\phi,\omega^a\wedge (\nabla_{V_a}H^{1,0})\lrcorner\,\phi\rangle
&=\sum_{a=1}^n g_Q(\nabla_{\bar V_a}H^{0,1}, V_a) |V_a\lrcorner\,\phi|^2.
\end{align}
Now, we put $C =\min_{a,x} |V_a\lrcorner\,\phi|^2\ge 0$. Then from the equations above, we have
\begin{align}
\sum_{a=1}^n&\{\langle\omega^a\wedge (\nabla_{V_a}H^{1,0})\lrcorner\,\phi,\phi\rangle+\langle\phi,\omega^a\wedge (\nabla_{V_a}H^{1,0})\lrcorner\,\phi\rangle\}\notag\\
&\geq C\sum_{a=1}^n\{g_Q(\nabla_{V_a}H^{1,0},\bar V_a)+\overline {g_Q(\nabla_{V_a}H^{1,0},\bar V_a)}\}\notag\\
&\geq C\ {\rm div}_\nabla(\kappa_B^\sharp).
\end{align}
In the last inequality, we used  $\sum_{a=1}^n g_Q(\nabla_{V_a}H^{1,0},\bar V_a) ={\rm div}_\nabla(H^{1,0})$. 
Since $\delta_B\kappa_B = -{\rm div}_\nabla(\kappa_B^\sharp) + |\kappa_B|^2$ and $\delta_B\kappa_B=0$, we have
\begin{align}
\sum_{a=1}^n&\{\langle\omega^a\wedge (\nabla_{V_a}H^{1,0})\lrcorner\,\phi,\phi\rangle+\langle\phi,\omega^a\wedge (\nabla_{V_a}H^{1,0})\lrcorner\,\phi\rangle\}\geq C |\kappa_B|^2.
\end{align}
Let $\phi$ be a basic-harmonic form of type $(r,0)$. Since the transversal Ricci curvature is nonnegative, from 
Corollary~\ref{AnotherLaplaceCor} and the inequality above, we have 
\begin{align*}
\Delta_B |\phi|^2 \leq 0.
\end{align*}
Hence by the maximum principle, $\phi$ is parallel. Since $Ric^Q>0$ at some point, $\phi$ is zero. $\Box$

\begin{coro} \label{VanishingCor}
Let $\mathcal F$ be as in Theorem~\ref{VanishingThm}. If the transversal Ricci curvature is nonnegative and positive at some point, then for $r>0$ there are no non-zero basic holomorphic forms of type $(r,0)$.
\end{coro}
{\bf Proof.} From \ref{FVanishingThm}, if the transversal Ricci curvature is nonnegative and positive at some point, then $\mathcal H_B^1=\{0\}$. Hence for some basic function $f$, $\kappa=df$, which implies $\Delta_Bf=0$. Hence by the maximum principle, $h$ is constant. That is, $\kappa=0$.    It follows from Corollary~\ref{basicholo(r,0)}. $\Box$

A theorem similar to the corollary above, with a stronger hypotheses on the transverse Ricci curvature, was proved in 
\cite[Theorem 3.4]{HV}.

\end{document}